\baselineskip=20pt

\font\twelvebf=cmbx12

\bigskip
\bigskip

\def\Z{{\bf Z}}
\def\C{{\bf C}}
\def\R{{\bf R}}
\def\Ch{{\C[[h]]}}
\def\N{{\bf N}}

\def\e{\eqno}
\def\l{\ldots}
\def\n{\noindent}

\def\t{\theta}

\def\q{{\bar q}}

\def\b{\bigskip}

\def\>{\rangle}
\def\<{\langle}
\def\v{\varphi}
\def\[{[\![}
\def\]{]\!]}
\def\Uq{U_q[gl(n/m)]}
\def\W{W(n-1/m)}
\def\A{{\tilde A}}

\hfill {{\twelvebf Preprint UR-1525 $\;\;$ 4/98}

\vskip 2cm

\noindent
{\twelvebf A Dyson realization and a Holstein-Primakoff realization 
for the

\noindent
quantum superalgebra U$_q$[gl(n/m)]  }

\vskip 32pt
\noindent
{\bf Tchavdar D. Palev}\footnote*{Permanent address: 
Institute for Nuclear Research and Nuclear Energy, 1784 Sofia, 
Bulgaria; $e$-mail: tpalev@inrne.acad.bg}

\noindent
Department of Physics and Astronomy, University of Rochester,
Rochester, New York 14627
\vskip 24pt

\noindent 
{\bf Abstract.} The Holstein-Primakoff and the Dyson realizations of
the Lie superalgebra $gl(n/m)$ are generalized to the class of the
quantum superalgebras $U_q[gl(n/m)]$ for any $n$ and $m$. It is shown
how the elements of $\Uq$ can be expressed via $n-1$ pairs of Bose
creation and annihilation operators and $m$ pairs of Fermi creation
and annihilation operators.

\vfill\eject 
\leftskip 0pt

\n 
Recently an analogue of the Dyson (D) and of the Holstein-Primakoff
(HP) realization for the superalgebras $gl(n/m)$ of any rank was given
[1]. In the present paper the results are extended to the quantum
superalgebras $U_q[gl(n/m)]$. The elements of $U_q[gl(n/m)]$ are
expressed via $n-1$ pairs of Bose creation and annihilation operators
(CAOs) and of $m$ pairs of Fermi CAOs. In the case $m=0$ the results
reduce to those announced in [2], namely to D and HP realizations of
$U_q[gl(n)]$ in terms of $n-1$ pairs of only Bose operators.

Initially the D and the HP realizations were given for $sl(2)$ [3,
4]. The generalization for $gl(n)$ is due to Okubo [5]. The ``quantum
case'' was worked out also first for $U_q[sl(2)]$ [6] and $U_q[sl(3)]$
[7]. Very recently it was extended to $U_q[sl(n)]$ [2]. To the best of
our knowledge analogues of D and of HP realizations for quantum
superalgebras have not been published in the literature so far. The
available realizations are of Jordan-Schwinger type, requiring $n$
pairs of Bose CAOs and $m$ pairs of Fermi CAOs in case of
$U_q[gl(n/m)]$ [8].

The motivation in the present work stems from the various applications
of the Holstein-Primakoff and of the Dyson realizations in theoretical
physics. Beginning with [2] and [3] the HP and D realizations were
constantly used in condensed matter physics. Some early applications
can be found in the book of Kittel [9] (more recent results are
contained in [10]). For applications in nuclear physics see [11, 12]
and the references therein, but there are, certainly, several other
publications.  Once the $q$-analogues of D and of HP realizations for
$U_q[sl(2)]$ and $U_q[sl(3)]$ were established, they found also
immediate applications [13-18]. One can expect therefore that the
generalization of the results to an arbitrary rank $U_q[gl(n/m)]$
superalgebra may prove useful too.

To begin with we recall the definition of $U_q[gl(n/m)]$ in the sense
of Drinfeld [19], keeping close to the notation in [20]. Let $\Ch$ be
the complex algebra of formal power series in the indeterminate
$h$, $q=e^{h/2}\in \Ch$. Then $U_q[gl(n/m)]$ is a Hopf algebra, which
is a topologically free $\Ch$ module (complete in the $h-$adic
topology), with generators $h_{j}$, ($j=1,2,\ldots,r\equiv n+m$) and
$e_i$, $f_i$ ($i=1,2,\ldots,r-1$) subject to the following relations
(unless stated otherwise, the indices below run over all possible
values)~:

\n
The Cartan-Kac relations~:
$$
\eqalignno{
& [h_{i},e_j]=(\delta_{ij}-\delta_{i,j+1})e_j;& (1)\cr
& [h_{i},f_j]=(-\delta_{ij}+\delta_{i,j+1})f_j;& (2)\cr
& e_if_j-f_je_i=0, \quad if \quad i\ne j; &(3)\cr
& e_if_i-f_ie_i={{q^{h_i-h_{i+1}}-q^{-h_i+h_{i+1}}}\over{q-q^{-1}}},
\quad if \quad i\ne n;& (4)\cr
& e_nf_n+f_ne_n={{q^{h_n+h_{n+1}}-q^{-h_n-h_{n+1}}}\over{q-q^{-1}}}; 
& (5)\cr
}
$$
The Serre relations for the $e_i$ ($e$-Serre relations)~:
$$
\eqalignno{
& e_ie_j=e_je_i, \quad if \quad  |i-j| \neq 1; \quad e_n^2=0; & (6)\cr
& e_i^2e_{i+1}-(q+q^{-1})e_ie_{i+1}e_i+e_{i+1}e_i^2=0, &\cr
& {\rm for} \quad i\in \{1,\ldots,n-1\}\cup\{n+1,\ldots,n+m-
2\};&(7)\cr
& e_{i+1}^2e_i-(q+q^{-1})e_{i+1}e_ie_{i+1}+e_ie_{i+1}^2=0,&\cr
&{\rm for }\quad i\in\{1,\ldots,n-2\}\cup\{n,\ldots,n+m-2\};& (8)\cr
& e_ne_{n-1}e_ne_{n+1}+e_{n-1}e_ne_{n+1}e_n+e_ne_{n+1}e_ne_{n-1}&\cr
& +e_{n+1}e_ne_{n-1}e_n-(q+q^{-1})e_ne_{n-1}e_{n+1}e_n=0;& (9)\cr
}
$$
The relations obtained from (6)-(9) by replacing
every $e_i$ by $f_i$ ($f$-Serre relations).

Let
$$
\t_i=\cases {{\bar 0}, & if $\; i<n$;\cr 
               {\bar 1}, & if $\; i\ge n$.\cr }\e(10)
$$
Then
$$
deg(h_i)={\bar 0}, \quad
deg(e_j)=deg(f_j)=\t_{j-1}+\t_j, \e(11)
$$
i.e. the generators $e_n$ and $f_n$ are odd and all other
generators are even.

We do not write the other Hopf algebra maps $(\Delta,\; \varepsilon,
S)$ (see [20]), since we will not use them. They are certainly also a
part of the definition.

The Dyson and the Holstein-Primakoff realizations are different
embeddings of 
$U_q[gl(n/m)]$ into the algebra
$W(n-1/m)$ of $n-1$ pairs of Bose CAOs and $m$
pairs of Fermi CAOs. The precise definition of $W(n-1/m)$
is the following. Let $A_{1}^\pm,\l, A_{n+m-1}^\pm$ 
be $\Z_2-$graded
indeterminates:
$$
deg(A_{i}^\pm)=\t_i.\eqno(12)
$$
Then $W(n-1/m)$ is a topologically free $\Ch$ module and 
an associative unital superalgebra with generators
$A_{1}^\pm,\l, A_{n+m-1}^\pm$ subject to the relations
$$
\[A_i^-,A_j^+\]=\delta_{ij},\quad
\[A_i^+,A_j^+\]=\[A_i^-,A_j^-\]=0. \eqno(13)
$$
Here and throughout
$$
\[x,y\]=xy-(-1)^{deg(x)deg(y)}yx,\quad 
\[x,y\]_q=xy-(-1)^{deg(x)deg(y)}qyx \eqno(14)
$$
for any two homogeneous elements $x$ and $y$.  With respect to the
supercommutator $\[x,y\]$ $W(n-1/m)$ is also a Lie superalgebra.

From (13) one concludes that $A_1^\pm,\ldots A_{n-1}^\pm $ are Bose
CAOs, which are even variables; $A_{n}^\pm,\ldots A_{n+m-1}^\pm $
are Fermi CAOs, which are odd. The Bose operators commute with the
Fermi operators.

In the physical applications it is often more convenient
to consider $h$ and $q$ as complex numbers, $h, q\in \C$. Then all our
considerations remain true provided $q$ is not a root of 1. The
replacement of $q\in \Ch$ with a number corresponds to a factorization
of $\Uq$ and $\W$ with respect to the ideals generated by the
relation $q=number$. The factor-algebras $\Uq$ and $\W$ are
complex associative algebras. However the completion in the $h$-adic
topology leaves a relevant trace: after the factorization the
elements of $\Uq$ and of $\W$ are not simply polynomials of
their generators. In particular the functions of the CAOs, which
appear in the D and in the HP realizations (see  bellow)
are well defined as elements from $\W$.

Now we are ready to state our main results. Let
$$
\q=q^{-1},\quad
[x]={{q^x-\q^x}\over{q-\q}}, \quad
N_i=A_i^+A_i^-, \quad N=N_1+\l+N_{n+m-1}.\e(15)
$$

\smallskip\noindent
{\it Proposition 1 (Dyson realization)}. The linear map 
$\varphi:\Uq \rightarrow \W$, defined on the generators as 
$$
\eqalign{
& \v(h_1)=p-N, \quad \v(h_i)=N_{i-1}, 
  \quad i=2,\l,n+m\equiv r,  \cr
& \v(e_1)={[N_1+1]\over{N_1+1}}[p-N]A_1^-,\quad
 \v(e_i)={[N_i+1]\over{N_i+1}}A_{i-1}^+A_i^-,
\quad i=2,\ldots,n-1,\cr
& \v(e_i)=A_{i-1}^+A_i^-,\quad i=n,\l,n+m-1,   \cr
& \v(f_1)=A_1^+,\quad
\v(f_i)={[N_{i-1}+1]\over{N_{i-1}+1}}A_i^+A_{i-1}^-,
\quad i=2,\ldots,n,  \cr
& \v(f_i)=A_i^+A_{i-1}^-,\quad i=n+1,\l,n+m-1, \cr
}\e(16)
$$
is a homomorphism of $\Uq$ into $\W$ for any $p\in \C$.

The proof is straightforward. One has to verify that Eqs.
(1)-(9) with $\v(h_i)$, $\v(e_i)$, $\v(f_i)$ substituted
for $h_i$, $e_i$ and $f_i$, respectively, hold. 
In the intermediate computations the relations
$$
[N_i,A_i^\pm]=\pm A_i^\pm,\quad [N,A_i^+A_j^-]=0.\e(17)
$$
$$
F(N_i)A_i^+=A_i^+F(N_i+1),\quad
F(N_i)A_i^-=A_i^+F(N_i-1),\e(18)
$$
are repeatedly used. The verification of Eqs. (4), (5) and (7)-(9) is
based also on the identity\hfill\break
$
[x+1]-(q+\q)[x]+[x-1]=0. 
$
Other relations, used only in the proof of (5)
are
$$
q^{N_n}=1-{N_n}+q{N_n}, \quad \q^{N_n}=1-{N_n}+\q {N_n}. \e(19)
$$

Similar as for $gl(n-1/m)$ [1],
the Dyson realization defines an infinite-dimensional
representation of $\Uq$ (for $n>0$) in the Fock space
${\cal F}(n-1/m)$ with orthonormed basis
$$
|l)\equiv |l_1,\ldots,l_{r-1})=
{(A_1^+)^{l_1}\ldots (A_{r-1}^+)^{l_{r-1}}\over
\sqrt{l_1!\ldots l_{r-1}!}}|0\>, \e(20)
$$
where
$
l_1,\ldots, l_{n-1}\in \Z_+\equiv \{0,1,2,\l\};
\;\; l_{n},\ldots, l_{r-1}\in \Z_2 \equiv  \{{\bar 0},{\bar 1}\}. 
$

If $p$ is a positive integer, $p\in \N$, the representation is
indecomposible: the subspace 
$$
{\cal F}_1(p;n-1,m)=
lin.env.\{\;|l_1,\ldots,l_{r-1}\;\>|\;l_1+\l+l_{r-1}> p\;\}\e(21)
$$
is an invariant subspace, whereas its orthogonal
compliment
$$
{\cal F}_0(p;n-1,m)=
lin.env.\{\;|l_1,\ldots,l_{r-1}\>\;|\;l_1+\l+l_{r-1}\le p\;\}\e(22)
$$
is not an invariant subspace. If $p\notin \N$, the representation is
irreducible. In all cases however, and this is the disadvantage of the
D realization, the representation of $\Uq$ in ${\cal F}(n-1,m)$ is not
unitarizable with respect to the antilinear anti-involution 
$\omega:\Uq \rightarrow \Uq$, defined on the generators as
$$
\omega(h_i)=h_i,\quad \omega(e_i)=f_i.\e(23)
$$

In order to turn ${\cal F}_0(p;n-1,m)$ into an unitarizable $\Uq$
module we pass to introduce the HP realization. To this end set
$$
\<N_i+c\>=\Bigg({[N_i+c]\over{N_i+c}}\Bigg)^{ {1-\t_i}\over 2}.
\e(24)
$$

\smallskip\noindent
{\it Proposition 2 (Holstein-Primakoff realization)}. The linear map 
$\pi:
\Uq \rightarrow \W,
$ 
defined on the generators as:
$$
\eqalign{
&\pi(h_1)=p-N, \quad \pi(h_i)=N_{i-1}, \quad i=2,\l,n+m\equiv r,  \cr
&\pi(e_1)=\sqrt{[p-N]}\;\<N_1+1\>A_1^-,\;\quad
\pi(f_1)=\sqrt{[p-N+1]}\;\<N_1\>A_1^+,\cr
& \pi(e_i)=\<N_{i-1}\>\<N_{i}+1\>A_{i-1}^+A_i^-,
\;\; i=2,\l,r-1, \cr
& \pi(f_i)=\<N_{i-1}+1\>\<N_{i}\>A_i^+A_{i-1}^-,
\quad i=2,\l,r-1 \cr
}\e(25)
$$
is a homomorphism of $\Uq$ into $\W$.
If $p\in\N$,
then ${\cal F}_0(p;n-1,m)$ and ${\cal F}_1(p;n-1,m)$ are invariant
subspaces; ${\cal F}_0(p;n)$ carries a finite-dimensional
irreducible representation; it is unitarizable with
respect to the anti-involution (23) and the metric defined with the
orthonormed basis (20), provided $h\in \R$. The representations,
corresponding to different $p\in \N$ are inequivalent

We scip the proof.  The circumstance that ${\cal F}(n-1,m)$ is a
direct sum of its invariant subspaces ${\cal F}_0(p;n-1,m)$ and ${\cal
F}_1(p;n-1,m)$ is due to the the factor $\sqrt{[p-N]}$ in $\pi(e_1)$
and $\sqrt{[p-N+1]}$ in $\pi(f_1)$. If $h\in \R$, then ($(\;,\;)$ 
denotes
the scalar product)
$$
(\pi(h_i)|l\>,|l'\>)=(|l\>,\pi(h_i)|l'\>),\quad
(\pi(e_i)|l\>,|l'\>)=(|l\>,\pi(f_i)|l'\>) 
$$
for all $|l\>, |l'\>\in {\cal F}_0(p;n-1,m),\;\;i=1,\l,r-1$.
Hence the representation of $\Uq$ in  ${\cal F}_0(p;n-1,m)$
is unitarizable.

Let us say a few words about the place of the Fock 
representations  among all known representations.  Any highest weight
finite-dimensional irreps $\psi$ of $gl(n/m)$ or $\Uq$ is labeled by
its signature $\{m\}\equiv \{m_{1},m_{2},\ldots,m_{r}\}$, where each
$m_{i}$ is determined from $\psi(h_i)x_0=m_{i}x_0$. Here $x_0$ is the
highest weight vector. So far explicit expressions for all
(finite-dimensional) irreps are available only for $gl(n/m)$ with
$m=1$ [21, 22]. Each such representation can be deformed also to an
irreps of $U_q[gl(n/1)]$ [23].

In case of $gl(n/m)$ or $\Uq$ with $m\ne 1$ explicit constructions 
were carried out for the so called essentially typical representations
[24, 20]. A representation is essentially typical, if
$$
\{l_{1},l_{2},\l,l_{n}\}\cap 
\{l_{n+1},l_{n+1}+1,l_{n+1}+2,\l,l_{r}\}=\emptyset,\e(26)
$$
where $l_i=m_i-i+n+1$ for $1\le i\le n$ and 
$l_j=-m_j+j-n$ for $n+1\le j \le r$.

In case of  ${\cal F}_0(p;n-1,m)$ the highest weight vector
is the vacuum. Then $m_i=p\delta_{1i},\;i=1,\l,r$ and therefore
(26) is not fulfilled. Hence the Fock space representations
of $\Uq$ in  ${\cal F}_0(p;n-1,m), \; p\in \N$ describe 
finite-dimensional irreps in addition to those studied in [20].
As mentioned already, if $q$ is taken to be a number, it
should be not a root of 1.

In conclusion we note that the operators [25]
$$
\A_i^-=\<N_i+1\>A_i^-, \quad \A_i^+=\<N_i\>A_i^+, \quad
{\tilde N}_i=N_i,
\quad i=1,\l,n+m-1 \e(27)
$$
satisfy the relations
$$
\[\A_i^-,\A_j^+\]_q=\delta_{ij} q^{-{\tilde N}_i},\quad
[{\tilde N}_i,\A_j^\pm]=\pm \delta_{ij}\A_j^\pm,\quad
\[\A_i^\pm,\A_k^\pm\]=[{\tilde N}_i,{\tilde N}_k]=0.\quad
i\ne k. \e(28)
$$ 
Therefore $\A_i^\pm,\;i=1,\l,n-1$, give a representation
of the algebra of the deformed Bose operators [26-28],
whereas $\A_i^\pm,\;i=n,\l,n+m-1$, yield a representation
of the deformed Fermi operators [8]. In terms of the
deformed operators Eqs. (25) read:
$$
\eqalignno{
&\pi(h_1)=p-{\tilde N}, \quad 
\pi(e_1)=\sqrt{[p-{\tilde N}]}\;\A_1^-,\quad
\pi(f_1)=\sqrt{[p-{\tilde N}+1]}\;\A_1^+,& (29a)\cr
& \pi(h_i)={\tilde N}_{i-1},\quad 
  \pi(e_i)=\A_{i-1}^+\A_i^-,\quad
  \pi(f_i)=\A_i^+\A_{i-1}^-,
\quad i\ne 1. & (29b)\cr
}
$$
The equations (29) could be called a $q-$deformed analogue
of the Holstein-Primakoff realization for $gl(n/m)$ [1],
whereas only Eqs. (29b) correspond to a $q$-deformed Jordan-Schwinger
realization of $U_q[gl(n-1/m)]$ [

\b\b\n
The author is thankful to Prof. S. Okubo for the kind invitation to
conduct the research under the Fulbright Program in the Department of
Physics and Astronomy, University of Rochester. This work was
supported by the Fulbright Program, Grant No 21857.

\vfill\eject

\n
{\bf References}

\b
\+ [1] & Palev T D 1997 {\it J. Phys. A: Math. Gen.} {\bf 30} 
        8273 \cr  

\+ [2] & Palev T D 1998 math.QA/9804017 \cr

\+ [3] & Holstein T and Primakoff H 1949 
         {\it Phys. Rev.} {\bf 58}  1098 \cr            

\+ [4] & Dyson F J 1956 {\it Phys. Rev.} {\bf 102} 1217 \cr

\+ [5] & Okubo S 1975 {\it J. Math. Phys.} {\bf 16}  528 \cr

\+ [6] & Chaichian M, Ellinas D and Kulish P P
         1990 {\it Phys. Rev. Lett.} {\bf 65}  980  \cr

\+ [7] & da-Providencia J 1993
         {\it J. Phys.  A: Math. Gen.} {\bf 26} 5845 \cr

\+ [8] & Floreanini R, Spiridonov V P and Vinet L 1990
        {\it Comm. Math. Phys.} {\bf 137} (1991) 149\cr  

\+ [9] & Kittel C 1963 {\it Quantum Theory of Solids} (Willey, 
          New York) \cr

\+ [10] &$\;$ Caspers W J 1989 {\it Spin Systems} (World Sci. Pub. 
Co., Inc.,
          Teanek) \cr

\+ [11] &$\;$ Klein A and Marshalek E R 1991
         {\it Rev. Mod. Phys.} {\bf 63} 375 \cr

\+ [12] &$\;$ Ring P and Schuck P  {\it The Nuclear Mani-Body 
Problem},\cr
\+      &    (Springer-Verlag, New York, Heidelberg, Berlin) \cr

\+ [13] &$\;$ Quesne C 1991 {\it Phys. Lett. A} {\bf 153} 
              303 \cr

\+ [14] &$\;$ Chakrabarti R and Jagannathan R 1991
         {\it J. Phys.  A: Math. Gen.} 
         {\bf 24}   L711 \cr

\+ [15] &$\;$ Katriel J and Solomon A I 1991
         {\it J. Phys. \ A: Math. Gen.} {\bf 24} 2093 \cr

\+ [16] &$\;$ Yu Z R 1991
         {\it J. Phys.  A: Math. Gen.} {\bf 24} L1321 \cr

\+ [17] &$\;$ Kundu  A and Basu Mallich B 1991 {\it Phys. Lett. A} 
          {\bf 156} 175  \cr

\+ [18] &$\;$ Pan F 1991  {\it Chin. Phys. Lett.} {\bf 8}  56 \cr

\+ [19] &$\;$ Drinfeld V 1986 {\it Quantum Groupd (Proc. Int. 
         Congress of Mathematics (Berkeley, 1986))},\cr
\+      & Ed.Gleasom A M (Providence, RI: American Physical
          Society) p 798  \cr

\+ [20] &$\;$ Palev T D, Stoilova N I and Van der Jeugt J 1994
        {\it Comm. Math. Phys.} {\bf 166} 367\cr 

\+ [21] &$\;$ Palev T D 1987 {\it Funct. Anal. Appl.} {\bf 21}  245  
            (English translation)\cr

\+ [22] &$\;$ Palev T D 1989 {\it Journ. Math. Phys.} {\bf 30} 1433\cr

\+ [23] &$\;$ Palev T D and Tolstoy V N 1991
            {\it Comm. Math. Phys.}  {\bf 141} 549\cr

\+ [24] &$\;$Palev T D 1989 {\it Funct. Anal. Appl.} {\bf 23} 
          141  (English translation).\cr

\+ [25] &$\;$ Polychronakos A P 1990 {\it Mod. Phys. Lett.}
            {\bf 5} 2325 \cr

\+ [26] &$\;$ Macfarlane A J 1989  {\it J.\ Phys. A : Math.\ Gen.}
          {\bf 22} 4581  \cr

\+ [27] &$\;$ Biedenharn L C 1989 {\it J.\ Phys. A: Math.\ Gen.} 
          {\bf 22}  L873 \cr

\+ [28] &$\;$ Sun C P and Fu H C 1989 {\it J.\ Phys. A: Math.\ Gen.} 
         {\bf 22}  L983  \cr

\end